\documentclass[11pt]{article}

\makeatletter
\@addtoreset{figure}{section}
\def\thefigure{\thesection.\@arabic\c@figure}
\def\fps@figure{h, t}
\@addtoreset{table}{section}
\def\thetable{\thesection.\@arabic\c@table}
\def\fps@table{h, t}
\@addtoreset{equation}{section}

\makeatother

\newcommand{\R}{\hbox{\bf R}}
\newcommand{\C}{\hbox{\bf C}}

\newtheorem{proposition}{Proposition}
\newtheorem{theorem}{Theorem}

\newtheorem{lemma}{Lemma}
\newtheorem{definition}{Definition}

\newtheorem{remark}{Remark}

\newtheorem{conjecture}{Conjecture}
\newtheorem{corollary}{Corollary}

\newtheorem{question}{Question}
\newtheorem{prescription}{Prescription}

\begin{document}

{\small
\title{Vortices on closed surfaces}
\author{ Stefanella Boatto\\
Departamento de Matem\'atica Aplicada da UFRJ \\ C.P. 68530, Cidade Universit\'aria \\ Rio de Janeiro, RJ, 21945-970, Brazil  (boatto@impa.br) \\
Jair Koiller\thanks{Supported in part by a grant from  Nehama gest\~ao de ativos.} \\
N\'ucleo de Matem\'atica Aplicada, Funda\c c\~ao Getulio Vargas \\ Praia de Botafogo 190 \\ Rio de Janeiro, RJ, 22250-040, Brazil  (jkoiller@fgv.br) }
}
\date{}
\maketitle
\pagestyle{headings}
{


{\scriptsize
\tableofcontents
}

\newpage

 \small{

\begin{center}
{\bf Abstract}
\end{center}

\begin{center}
{ \footnotesize
}
\end{center}

\noindent We consider $N$ point vortices $s_j$  of strengths $\kappa_j$ moving  on  a closed
(compact, boundaryless, orientable) surface  $S$ with riemannian metric $g$.  As far as we know, only the sphere or surfaces of revolution, the latter
qualitatively, have been treated in the available literature. The aim of this note is to present an intrinsic geometric formulation for the
general case.
Since the pioneer works of Bogomolov \cite{Bogomolov} and Kimura/Okamoto \cite{KiO} on
the sphere $S^2$, it is known that stream function produced by a unit point vortex at $s_o \in S$ on a background uniform counter vorticity field is given by Green's function  $G_g(s,s_o)$ of the Laplace-Beltrami operator
$\Delta_g = {\rm div}g \circ {\rm grad}_g$.
It  behaves as $  \log d(s,s_o)/2\pi$ near $s_o$.  Desingularizing  the stream function $G_g$  is therefore mainstream mathematics (no pun intended).
An energy core argument shows that a single   vortex $s_o$ drifts according to the Hamiltonian system $(\Omega_g,R_g)$,
where $\Omega_g$ is the area form  of $g$ and $R_g$ is Robins's function $R_g(s_o) = \lim_{s \rightarrow s_o} \, G_g(s,s_o) - \log d(s,s_o)/2\pi$ .
 If $S$ has genus zero,
it it known that $R_g = \Delta_g^{-1} K + {\rm trace} \Delta_g^{-1}/A(S)$.
The collective Hamiltonian keeps
the same form of C.C.Lin's now classical paper \cite{Lin}, with
$$\Omega_{{\rm collective}} = \sum_{j=1}^N \, \kappa_J \Omega(s_j)\,\,\,,\,\,\,
 H = \sum_{1 \leq i < j \leq N} \kappa_i \kappa_j G_g(s_i,s_j) + \sum_{\ell=1}^N   \frac{1}{2} \kappa_{\ell}^2  R_g(s_{\ell}) \,\,.$$
  For Jordan domains $D \subset S$ the structure is the same,  using the appropriate hydrodynamical Green function.   The extension for vortices with mass is also immediate.                 Under conformal changes of metrics $\tilde{g} = h^2 g$
the symplectic form changes accordingly to  $ \Omega_{\tilde{g}} = h^2 \Omega_g$ while the new Hamiltonian is given by
$$
\tilde{H} = H - \frac{1}{4\pi} \sum_{\ell=1}^N \kappa_{\ell}^2 \, \log(h(s_{\ell})) - \frac{\kappa}{\tilde{A}(S)} \sum_{\ell=1}^N \kappa_{\ell} \Delta^{-1} h^2(s_{\ell})\,\,\,\,,\,\,\,  \kappa = \sum_{\ell=1}^N \kappa_{\ell} \,\,.
$$
  The presence of the total vorticity in the
  conformal change formula reflects the fact that when the sum of the vorticities vanish, the collective vortex stream function
$\psi(s; s_1,...,s_N)$ is independent of the conformal metric $\tilde{g} = \exp(2\phi) g$. In this case the $\tilde{g}$-regularization of $\psi$
at any of the vortices simplifies:   one just subtracts off $\frac{1}{2\pi} \phi(s_j)$ from the $g$-regularized
stream function at $s_j$. In particular, when $S$ is conformal to the standard sphere, making
an artificial puncture at any given $s^* \in S$ allows to easily write vortex motions for any metric on the sphere.  We give a simple proof of Kimura's conjecture that a dipole describes   geodesic motion. Searching for integrable vortex pairs systems on Liouville surfaces is in order. The vortex pair system on a triaxial ellipsoid extends Jacobi's
geodesics. Is it Arnold-Liouville integrable?   Not in our wildest dreams is another possibility:  that quantizing a vortex system could relate  with a million dollars worth question, but we took courage - nerve is more like it - to also present it.
}



\newpage

\footnotesize{

\section{Introduction }

The purpose if this note is to formulate vortex motions in   intrinsic  fashion. The (now classical) results
of C.C.Lin (\cite{Lin}, 1941) are generalized and reinterpreted.
Complex coordinates $z \in \C $ will not be used (well, almost).  The key observation  is that the core energy desingularization procedure (\cite{Gustafsson}, 1997)   geometrizes in a nice way,  matching perfectly with  Green's function behavior.

\subsection{Vortex history on a capsule}


   Helmholtz' ``Wirbel'' paper, short for  ({\it \"Uber integrale der hydrodynamischen gleichungen
   welche den Wirbelbewegungen entsprechen}, \cite{Helmholtz}, 1858) came out a century$+$ after Euler's {\it Principes g\'en\'eraux
du mouvement des fluides} (\cite{Euler}, 1757), where  the  mathematical definition of
vorticity first appeared\footnote{Circa 340 B.D. Aristoteles  tried to explain
  typhoons in his {\it Metereologica}.  The importance of  vorticity in microfluid motion was clearly recognized
by  Leonardo da Vinci.  Decartes envisaged vortices as powerhouses for planetary
movement. W. and J.J. Thomson used vortices for an  atomic theory, discarded after Rutherford's 1912 experiments. Some say string theory  is a souped up version of the primeval vortex soup.}.
 Eighteen years later Kirchhoff presented the
equations for point vortices in the plane ({\it Vorlesungen}, \cite{Kirchhoff}, chapter XX,  p. 259, eq.(14)):
$$
m_j \frac{dx_1}{dt} =  \frac{\partial P}{\partial y_j} \,\, , \,\,m_j  \frac{\partial y_1}{\partial t} =  - \frac{\partial P}{\partial  x_j} \,\,\,\,\,\,\,\,,\,\,\,\,
P =  \sum \frac{1}{\pi} \, m_i m_j \log  \rho_{ij} \,
$$
(nowadays  vortex strengths are  denoted by  $\kappa$ or $\Gamma/2\pi$ instead of $m$,     $1/\pi$  is omitted, and     $H$
replaces $P$).

This is a Hamiltonian system $(\Omega, P)$, where the symplectic form  is the combination of planar areas weighted by the vorticities,
$  
\Omega_{system} =  \sum_j \,  \kappa_j dx_j \wedge dy_j \,\,.$    
Only in 1941 the correct Hamiltonian was found by C.C.Lin, \cite{Lin}, for  the case of
 arbitrary planar domains  $D\, \subset \C$.   See \cite{Newton} and \cite{BM} for treatises and \cite{Aref}, \cite{Aref1}, \cite{BoattoCrowdy} for recent reviews. 

\subsection{Vortices on curved surfaces}

Vortex dynamics on the sphere and hyperbolic plane started with Bogomolov \cite{Bogomolov}
and   Kimura and Okamoto \cite{Kimura}.   Basically, since the late 1970's up to now  only the constant curvature
case has been considered.   Follows a very incomplete list 
(we apologize for omissions): Kimura and Okamoto (\cite{KiO}, 1987),   Pekarski and Marsden (\cite{PeMa}, 1998), Kidambi and Newton (\cite{KiNa}, 1998, \cite{KiN},2000),  Lim et al. (\cite{Lim} 2001),
  Cabral and Boatto 
 (\cite{CaBo}, 2003), Borisov and al. (\cite{BoPa}, \cite{BoLe}, 1998,  \cite{BoKi}, 2000),
Laurent-Polz  (\cite{Laurent}, 2005), Tronin (\cite{Tronin}, 2006). The case of surfaces of revolution has just been touched
recently, in qualitative terms (\cite{Souliere}, 2002) or perturbatively (\cite{Castilho}, 2007).

 How about the general case?   Vortex equations
in the case of domains conformal to a  planar domain
were obtained by Hally \cite{Hally}.   The key aspect is how to desingularize   the {\it velocity} field at a vortex. Using isothermic coordinates, and after some   heuristics  (complicated for mathematicians)  
he found    additional   ``self-terms'' attributed to the
 {\it variations in  the surface curvature} near each  vortex.
For a  surface homeomorphic to
the z-plane $\C$, with metric $g = h^2(z,\overline{z}) dz \otimes d\overline{z}$,  Hally's equations are
\begin{equation}  \label{Hallyequations}
h^2(z_j,\overline{z}_j) \dot{\overline{z}}_j = i  \sum_{\ell,\ell \neq j}\, \frac{\kappa_\ell}{z_\ell-z_j}  \,\,+ \,\, i \kappa_j
\frac{\partial}{\partial z_j} \log (h(z_j,\overline{z}_j))\,\,\,,\,\,\,\, (1 \leq j \leq N) \,\,.
\end{equation}
A single unit vortex     obeys the ``self-motion equation''
  \begin{equation} \label{selfmotion}
   \dot{\overline{z}} =   i   h^{-2}(z,\overline{z})
\frac{\partial}{\partial z } \log (h(z,\overline{z})) \,\,.
   \end{equation}
    We provide corrections to \cite{Hally} for the case of nonzero total vorticity.  Our results are valid for closed surfaces of arbitrary genus and Jordan domains on  them.



 \section{Laplace-Beltrami operator and 2-dimensional hydrodynamics}

\subsection{Complex and symplectic geometry} \label{nonsense}

Hereafter  $S$ will always denote a two dimensional compact boundaryless orientable manifold (hereafter called closed),  endowed with a riemannian metric $g = \langle , \rangle$, and $D \subset S$ a Jordan domain (a region with compact closure on a whose boundary consists of a finite number of connected components each diffeomorphic to a circle).   The  underlying    {\it Riemann surface structure} for $S$ is given by the atlas ${\cal A}$  of Laplace-Beltrami
isothermic coordinates.  {\it Complex geometry }
keeps the notion of angle between tangent vectors $u,v$  (thus the rotation operator $J$ by $90$ degrees of tangent vectors is well defined), but neglect lenghts $ |u| $ and  areas
$\Omega(u,v)$.  {\it Symplectic geometry},  on the other hand,  keeps the area
form  but neglects the angles.
 Merging different symplectic forms with the same complex geometry produces
conformal riemannian metrics.  In 2-dimensions all   metrics  are automatically Kahler.

Functions (usually  velocity potentials   or stream functions) produce vectorfields under the
 gradient  and  symplectic gradient operators ${\rm grad}$ (usually denoted $\nabla$), $\, {\rm sgrad}: = J  \circ {\rm grad}$ (also denoted
  $ \nabla^{\perp}$),
\begin{equation} d\phi (\bullet) = \langle {\rm grad \phi} , \bullet \rangle = \Omega({\rm sgrad}(\phi), \bullet) \,\,.
\end{equation}
There is a
 ``baby Hodge theory'' in two dimensions: one identifies functions with two-forms via  $f \leftrightarrow  f \Omega$ and
 vectorfields with 1-forms via $ v \leftrightarrow \omega = \langle v , \bullet \rangle $.
 Let $C$ be a closed curve enclosing a domain homeomorphic to a disk $D \subset S$. Orient $C$ such that the frame
 $n,t$ is positive, where $t$ is the unit tangent vector and $n$ the exterior normal. Let $v$ be
 a vectorfield. The familiar formula
 \begin{equation}  \oint_C \langle v , n \rangle \, ds = \int\int_D  {\rm div}(v) \Omega
 \end{equation}
 makes perfect sense.  The coordinate free definition of the {\it divergence operator},
 sending vectorfields to functions, can be given as follows (applying Stokes' theorem and looking at the integrands):
\begin{equation}
d \langle J(v) , \bullet \rangle = {\rm div}(v) \Omega \,\,.
\end{equation}
{\it Incompressible hydrodynamics} means: ${\rm div} (v) = 0$  except for
 sources or sinks. On a compact surface $S$ their divergences must add to zero.
 In {\it potential} flow (and always incompressible), $ v = - {\rm grad} \phi$  for a velocity potential $\phi$ which in general is not  globally defined (it is a multivalued function in the Rieman surface sense).

 Define the Laplace operator by $ \Delta =  {\rm div } \circ {\rm grad} $.  It is a negative definite self-adjoint operator with respect to
$\langle f , g \rangle = \int_M  f g \,\Omega$. 
Incompressibility
implies  $\Delta \phi = 0 $ (up to singularities).
 Thus the velocity potential $\phi(s), s \in S$ is harmonic up to singularities, and  complex function theory of $S$ (viewed as a Riemann surface) comes to fore. The conjugate harmonic function $\psi$ is called the {\it stream} function.
Combined they form the {\it velocity potential}  $F = \phi + i \psi$ (in planar hydrodynamics, $u+iv = - \overline{F'(z)}$).   Clearly  $v$ is tangent to the level lines of $\psi$, and can
be described symplectically from $\psi$ via
\begin{equation} d\psi  = \Omega( v , \bullet ) \,\,\, \Leftrightarrow \,\,\, v = {\rm sgrad} (\psi) = - {\rm grad}(\phi) \,\,.
\end{equation}
The {\it circulation} $ \oint_C \, \langle v , t \rangle ds
$ of $v$ around a closed curve  $C = \partial D$
is an object of fundamental importance. Using Stokes' theorem it can be transformed into a double integral,
$$ \oint_C \, \langle v , t \rangle ds = \oint_C \, \langle J(v) , J(t)  \rangle ds = \oint_C \, \langle J(v) , -n  \rangle ds = - \int\int_D {\rm div}(Jv) \Omega
$$
But $J(v) = J (J {\rm grad} \psi) = - {\rm grad} \psi$. Hence
\begin{equation}  \oint_C \, \langle v , t \rangle ds =  \int\int_D  \Delta \psi \Omega
\end{equation}
A coordinate free definition of the surface {\it rotational} or {\it vorticity} $ \,\,\,\omega = {\rm rot v} \,\,\,$ is then
\begin{equation}   d \langle v , \bullet \rangle =    {\rm rot}(v) \Omega    \,\,.
\end{equation}

In passing, we mention that discrete versions of all these objects have been proved   useful in computer graphics and flow visualization, see e.g.  \cite{Dong}, \cite{Polthier}.

\subsection{Poisson's equation and  Green functions for the Laplace operator} 

The vorticity distribution $\omega$ forms ``the sinews and muscles of fluid motion'' (K\"uchemann, \cite{Kuchemann}, 1965).  No wonder   2-dimensional hydrodynamics is governed by Poisson's equation (inverting the Laplace operator). Perhaps the noblest
equation  of them all, it is the holy grail of every mathematician, pure or applied:
\begin{equation} \label{Poisson}
   \Delta \psi =  \omega \,\,, \,\,\,\,  \omega \in {\rm some} \, {\rm class} \, {\rm of}\, {\rm functions}\,\,.
\end{equation}
On a closed surface  the source term must average to zero:
\begin{equation}  \int\int_S  \omega \Omega = 0 \,\,.
\end{equation}
[ Proof. $\int\int_S\omega\Omega = \int\int_S\Delta\psi dS = \int\int_S \nabla\cdot \nabla \psi dS = \oint_{\partial S}
\langle \nabla\psi,n\rangle ds = 0$, since there is no boundary. Alternatively, take a small curve $C$ enclosing a disk $D$. The circulation of the velocity field can be computed as
the double integral of the vorticity on $D$ and minus that integral on $S-D$.]

 In two dimensions vorticity is a ``material property'',  $L_v \omega = 0$, and this translates into
\begin{equation}  \label{transport}
\frac{\partial \omega}{\partial t} + \langle v , {\rm grad} \omega \rangle = 0 \,\,,  v = {\rm sgrad} \psi \,\,.
\end{equation}

Choosing an appropriate class of vortex functions is an important mathematical problem.  In many concrete examples, the ``vorticity çountour map'' has sharply defined islands of high vorticity (positive or negative)
in  a sea of roughly constant (often zero) vorticity. In time, such regions can blend or develop  filaments. The evolution of
   the system (\ref{Poisson}, \ref{transport}) may reach a stand-off by the appearance of singularities.

  At any rate, 
 given $\kappa_j\,, 1 \leq j \leq N$, 
 take  geodesics disks $D_j$ of equal but very small area ${\rm Area}(Dj) =
\epsilon^2$
centered at $s_j$ and define
\begin{equation} \label{patch} \omega(s) =  \frac{\kappa_j}{\epsilon^2}    \,\,,\,\,\, s \in D_j
\end{equation}
and the suitable constant elsewhere.
One can prove that an initial state consisting of small nearly circular vortex cores remains stable.

 A {\it vortex pair} is the case $N=2$ with opposite vorticities. Their ubiquous presence in fluid flows   is  fascinating. In section \ref{Kimura} we give a simple proof of a conjecture by Kimura \cite{Kimura}:
a dipole
describes a geodesic on $S$ (zooming across the fluid very quicky).

A {\it monovortex} consists on a single vorticity island on  a  sea of distributed counter vorticity (uniformly, say) so that the average vorticity is
zero.
By linearity, any system of vortex patches can be decomposed on a weighted (positive and negative) sum of monovortices.
The following {\it gedanken} experiment may help. Consider a fluid at rest occupying the closed surface $S$. All of a sudden,
stir around a given point $s_o$, producing a vortex singularity there. What would be the resulting flow? Since there is
no distinguished point $s_*$ to be taken as a counter-vortex, the opposite vorticity should be distributed
as homogeneously as possible on the surface.  [ {\it Vorticity}  is an {\it area} related concept (via Hodge theory, vorticity is thought as a two-form), and
hence the area element of $S$ should be used to perform this homogenization.]  Note that a distinguished $s_*$ may exist if  $S$ has symmetries, for instance if $s_o$ is the south pole of a surface of
revolution; then $s_*$ is  the north-pole; those situations would be quite unstable, though\footnote{The lack of unicity
in some Euler flows is a known possibility; Joe Keller's ``tea pot effect'' is a  sometimes dramatic  example \cite{Joe}.}.

 Point vortex systems  are meant as finite dimensional ODE aproximations for   Euler's equations
 (geodesic flow on the group of area preserving transformations of $S$) or the equivalent PDEs (\ref{Poisson},\ref{transport}),  assuming concentrated vorticities
at  points $s_j(t), \, 1 \leq j \leq N. $
As $\epsilon \rightarrow 0$, how to obtain a suitable ODE approximation for the PDE system (\ref{Poisson}, \ref{transport})?   This is  done in two steps. Firstly, one needs the stream function for a marker particle on a flow generated by bound
   point vortices.  The second step is  a desingularization procedure in order to remove the influence
 of each vortex on itself.  Then we let them   ``dance'' under each other's influence and its own drift.

 The first step is achieved with Green functions od the Laplace-Beltrami operator.

\section{Point vortex equations}

\begin{definition}  \label{Green}
 $G = G_{(S,g)}$ is the Green function of $\Delta_g$,
\begin{equation} \Delta^{-1} \omega(s) = \int_M G(s,r) \omega(r) dA(r)\,\,,\,\,\,\, \Delta^{-1} \Delta \omega = \omega - \frac{1}{A} \int_S \omega dA \,\,.
\end{equation}
$G(s,s_o)$ is a smooth function outside the diagonal,   characterized by the following properties \cite{Okikiolu}:
\begin{equation} \label{characterization}
\Delta G(s,s_o)   =   - \frac{1}{{\rm Area(S)}}\,\,(s \neq s_o) \,\, , \,\, G(s,s_o) - \log d(s,s_o)/2\pi \,\,\, {\rm bounded}, \,\,\,\,\,\,
 \int_S G(p,q) \Omega(q)   =   0 \,\, , \,\,
 G(s,s_o)   =   G(s_o,s) .
\end{equation}
\end{definition}
Symmetry is to be expected from this formal calculation:
$$ G_{s_o}(s) = \int\int_S G_{s_o}(s_1) \Delta G_{s}(s_1) ds_1 = \int\int \nabla G_{s_o}\cdot \nabla G_{s}-
\int_{\partial S}G_{s_o}\partial_n G_{s}
= \int\int \nabla G_{s_o}\cdot \nabla G_{s} \,. $$

Hydrodynamically, $ \psi = G(s,s_o)$ fulfills all the requisites  to be the stream function of a unit monovortex at $s_o$.  In the   electrostatic interpretation,    $G(p;q)$ is the potential
at $p$ of a positive unit point charge at $q$, on a conductor with total charge zero (the background negative charge is uniformly
distributed).  A test particle $s$ with positive but negligible charge  on the field generated by a unit
charge will move according to
\begin{equation}  \label{charge}
\dot{s} = {\rm grad} G(s;s_o)\,\,,\,\, s \neq s_o \,\,.
\end{equation}
while a fluid particle (marker)  $s$ on the flow generated by an unit strength bound vortex $s_o$ will move according to
\begin{equation}  \label{marker}
\dot{s} = {\rm sgrad} G(s;s_o)\,\,,\,\, s \neq s_o \,\,.
\end{equation}

\subsection{Eppur si muove: vortex drift under its own stream}

Newton's third law precludes the action of the material particle at $s_o$ on itself, but   the rest of the universe $S - {s_o}$ conspires to
 impinge a  ``collective reaction'' on $s_o$ and makes it drift.
 In order to capture this motion, a
desingularization procedure is needed. To make that
long story short, a single vortex will move according to
\begin{prescription}  \label{desingularization}
\begin{equation}  \label{drift}
\dot{s}_o = {\rm sgrad} R(s_o)\,\,\,,\,\,\,
 R(s_o) = \lim_{s \rightarrow s_o}\,  G(s,s_o) - \frac{\log d(s,s_o)}{2\pi}
\end{equation}
where $d$ is the geodesic distance between points as measured with $g$.
\end{prescription}
The explanation for this rule is given in Proposition \ref{corenergy} below.
We were filled with  a religious feeling by the fact that the
regular part $R$ is  a notable object from the geometric function theory of the Laplace-Beltrami operator,
known as Robin's function of $g$ (\cite{Okikiolu2},
\cite{Okikiolu2}).      Vortex drift
belongs to mainstream mathematics!
\begin{remark} Jean Steiner \cite{Steiner} has shown that for any metric $g$ on $S^2$
\begin{equation} \label{Steiner}
R_g(s) - \frac{1}{2\pi}\, (\Delta_g^{-1} K_g)(s) = \frac{1}{A_g(S)} {\rm trace} \Delta_g^{-1}
\end{equation}
where $K$ is the Gaussian curvature.  For higher genus the RHS is no longer a constant and its flutuation does not
have a simple geometric interpretation.
\end{remark}
\begin{remark}
Marker's equation (\ref{marker}) becomes a 1 1/2 degrees of freedom system, by inserting  $s_o(t)$, a solution of
(\ref{drift}).  Chaotic marker equations are possible even for  one vortex systems.
\end{remark}
A basic observation  in this paper is  that the  core energy argument in \cite{Gustafsson} geometrizes nicely:
\begin{proposition} \label{corenergy}
The desingularization rule (\ref{drift}) encodes the core energy renormalization argument.
\end{proposition}
{\bf Proof.} Consider a small Gauss system of normal coordinates
around $s_o$.  The gradient ${\rm grad} d(s,s_o)$ is a unit vector along the geodesic rays. The symplectic gradient
${\rm sgrad}\, d(s,s_o)$ is obtained by composition with $J$ and is therefore  tangent to the geodesic circles.   Its flow leaves
the geodesic disks invariant, and the same is true for any function of the distance, in particular
$\psi_{{\rm sing}} = \log d(s,s_o)$.
For $\log d(s,s_o)$  the vectors   rotate with speed inversely proportional   to $d(s,s_o)$ on a geodesic circle. A simple back of the envelope reasoning shows that the kinetic energy confined in this
blob is logarithmically infinite. It can be removed (in physics jargon, renormalized) provided only
a small quantity of kinetic energy crosses the boundary of a geodesic disk in finite time. We invoke \cite{Gustafsson}:
energy diffusion ``can be neglected in the limit as the radius of the ball tends to zero'' ,
section 5). This  is true precisely due to the the fact that behavior of Greens's function $G(s,s_o) -  d(s,s_o)/2\pi$ is bounded.

\subsection{Collective vortex motions: Hamiltonian description}

Collective motions are obtained under
  \begin{prescription}
\begin{center}
 Each vortex   moves on
the stream generated by all the others  together with its own drift term.
\end{center}
\end{prescription}
Rewriting this prescription in  Hamiltonian form is   straightforward.  In view of   the   symmetry properties of Green's functions, C.C.Lin's \cite{Lin} arguments
geometrize directly.


\begin{theorem} \label{main}
The dynamics of $N$ point vortices $s_j$  of strenghts $\kappa_j$ moving on $S$ is a Hamiltonian system on $S \times ... \times S$ with symplectic form
\begin{equation}  \label{bigsymplectic}
\Omega_{collective} = \sum_{1 \leq j \leq N} \, \kappa_j \Omega(s_j) \,\,
\end{equation}
and Hamiltonian
\begin{equation} \label{bighamiltonian}
H = \sum_{i < j} \, \kappa_i \kappa_j  G(s_i,s_j) + \sum_{i=1}^N \, \frac{1}{2} \kappa_i^2 \, R(s_i) \,\,.
\end{equation}
where the self interaction terms are given by Robin's function
\begin{equation} \label{robin}
R(s) = \lim_{r \rightarrow s} \, G(r,s) - \frac{\log(d(r,s))}{2\pi} \,\,.
\end{equation}
\end{theorem}
\begin{remark} \mbox{}

\noindent 
i) The key aspect for deriving vortex equations is the geometric desingularization prescription (\ref{robin}).   As we
just observed, it can be solidly justified via the
core energy method, see section 5 of  Flucher and Gustafsson (\cite{Gustafsson}, 1997).  \\
ii) We omitted the possibility of collisions. See \cite{Lacomba1}, \cite{Lacomba2} for collisions and regularizations for 3 and 4 vortex motions on the plane.
\end{remark}

 There is an immediate extension to  vortices with mass  that
may be relevant due to   current interest in Bose-Einstein condensates,  see   \cite{KRO}, \cite{Borisovmass},
\cite{Ramodanov}).  These systems will  exhibit slow-fast Hamiltonian phenomena \cite{Neishtadt}.

\begin{corollary}
Let each $s_j$   have a mass $m_j$ besides its vorticity $\kappa_j$,    contributing with a kinetic energy
$  \frac{1}{2 m_j} \,  \langle  {p_s}_j ,  {p_s}_j  \rangle$,   where the bracket denotes the  the induced inner product in $T^*S$ via the
Legendre transform of $g$.
The  dynamics in  $T^*(S \times .... \times S)$ is described by the   Hamiltonian system
\begin{equation}    \label{2}
H =  \sum \,  \frac{1}{2 m_j} \,  \langle  {p_s}_j ,  {p_s}_j  \rangle \,\,  + \,\,  \sum_{i<j} \,  \kappa_i  \kappa_j \,  \hat{G}(s_i,s_j)\,\,\,\,\, \,\,\, , \,\,\,\,\,\,\,\,\,\,\,\,
 \Omega_{collective} =    \Omega_{can} +   \sum_j  \, \kappa_j \Omega(s_j)
\end{equation}
where   $\Omega_{can} =  \, ''  \sum \, dp_j \wedge ds_j  \,\,  '' $ is the canonical 2-form of  $T^* (S \times ... \times S)$.
\end{corollary}

We now address the following question. Let $\tilde{g} = h^2 g$  a conformal metric. How do  the vortex problems on $(S,g)$ and $(S,\tilde{g})$ compare?
Of course, the  symplectic form deforms as $ \Omega \rightarrow h^2\Omega$
in each coordinate.   Can we find the new  Hamiltonian in terms of the old one? We need transformation
formulas both  the Green and Robin functions.



\subsection{How Green and Robin functions change under conformal transformations}

Let  $\tilde{g} = h^2 g$ a conformal change of metric. Green's function  $G_{\tilde{g}}(s,s_o)$ for the
Laplace-Beltrami operator on a closed surface  is {\it not} conformally invariant. It {\it must} change,   because
a background uniform vorticity is an area dependent notion, except when it is zero. In fact, if the background vorticity is zero, the following ``naive conformal rule'' holds:
\begin{proposition} \label{naive}
 In a  zero vorticity background, in order to find the $\tilde{g}$-regularized stream function for a vortex, subtract $\log(h)/2\pi$  from its   $g$-regularized stream function,  $\tilde{g} = h^2 g$.
 \end{proposition}
{\bf Proof}.
$$ \psi_{{\rm reg}}(s_o) = \lim_{s \rightarrow s_o} \, \psi_{{\rm conformal}}(s) -   \log d(s,s_o)/2\pi \,\,.
$$
$$ \log(\tilde{d}(s,s_o)) = \log \left( d(s,s_o) \frac{\tilde{d}(s,s_o)}{d(s,s_o)} \right) \sim \log(d(s,s_o)) + \log h(s_o)  \,\,\,\,\, {\rm for}\,\, s \,\,{\rm near}\,\, s_o.
$$

However, {\it if, and this is a big if} the sum of vorticities does not vanish,  there is
a constant {\it nonzero} background counteracting vorticity in $S - \{s_1,...,s_n \}$.  Since what is constant is
a metric dependent notion,  this complicates
the conformal transformation rules.  Using  Green's function for the Laplace-Beltrami operator
on a closed surface   becames unavoidable.
The  correct transformation
 formulas are as follows  
 \begin{lemma} [\cite{Okikiolu}] \label{Okikiolu2}
 \begin{equation}
 \tilde{G}(s,s_o) - G(s,s_o) = - \frac{1}{\tilde{A}} \left( \Delta_g^{-1}h^2(s) +
 \Delta_g^{-1}h^2(s_o) \right) +
 \frac{1}{\tilde{A}^2}\, \int_S \,h^2 \Delta_g^{-1}h^2 \, \Omega
 \end{equation}
 \begin{equation}
\tilde{R}(s) = R(s) - \frac{1}{2\pi} h - \frac{2}{\tilde{A}} \Delta_g^{-1} h^2(s) + \frac{1}{\tilde{A}} \,\int_S \, h^2 \Delta_g^{-1} h^2 \Omega
\end{equation}
 \end{lemma}
A beautiful theory has been recently developed around Robins's function and   spectral invariants of conformal classes of metrics,
 see \cite{Steiner}. For our purposes the last term in both formulas can be dropped, as they just add  constants.  After some algebra juggling, one gets
 \begin{theorem} \label{confconf}
Under a conformal change of metric $\tilde{g}(s) = h^2 g(s)$, the  symplectic form deforms as $ \Omega \rightarrow h^2\Omega$
in each coordinate, and  the  Hamiltonian becomes
\begin{equation} \label{conformalrule}
\tilde{H} = H - \frac{1}{4\pi} \sum_{\ell=1}^N \kappa_{\ell}^2 \,\log(h(s_{\ell})) - \frac{\kappa}{\tilde{A}} \sum_{\ell=1}^N \kappa_{\ell} \Delta^{-1} h^2(s_{\ell})\,\,\,\,,\,\,\,  \kappa = \sum_{\ell=1}^N \kappa_{\ell} \,\,.
\end{equation}
 \end{theorem}
 \begin{remark}
 A map $F:(S,g) \rightarrow (\tilde{S},\tilde{g})$ is conformal when
 \begin{equation}
 \tilde{g}(dF_s \cdot v_s ,\,\, dF_s \cdot w_s ) = h^2(s) g(v_s,w_s) \,\,.
 \end{equation}
 The left hand side induces a metric $\tilde{g}_F(v_s,w_s)$ in $S$, hence  $\tilde{g}_F = h^2 g$ and
 Theorem \ref{confconf} applies. The vortex dynamics $(\tilde{\Omega},\tilde{H})$ on $S \times ... \times S$
 corresponds via $F$ to the  dynamics of vortices in $\tilde{S} \times .... \times \tilde{S}$.
 \end{remark}
\smallskip
For completeness, we present the
Proof of lemma \ref{Okikiolu2}},
following \cite{Okikiolu2}.  The trick is average over $\tilde{\Omega}$ twice using the  ``axioms'' in Definition \ref{Green}.   Consider the functions
$$G(s,s_o) - G(s,s_1) = E(s;s_o,s_1)\,\,,\,\, \tilde{G}(s,s_o) - \tilde{G}(s,s_1) = \tilde{E}(s;s_o,s_1) \,\,.$$
 Both $\tilde{E}$ and $E$ are harmonic up to (+) and (-) $\log$ singularities at $s_o$ and $s_1$
so they differ by a constant $c$. To find this constant we do the $\Omega(s) = h^2(s) \Omega(s)$ average of
 $$\tilde{E} - E =  (\tilde{G}(s,s_o) - \tilde{G}(s,s_1)) - (G(s,s_o) - G(s,s_1)) = c  \,\,.$$
 The first two terms drop out while the last two give
 $$  \Delta^{-1} h^2(s_1) - \Delta^{-1} h^2(s_o) = \, c \,\tilde{A} \,\,.
 $$
 Hence
 $$    (\tilde{G}(s,s_o) - \tilde{G}(s,s_1)) - (G(s,s_o) - G(s,s_1)) = \Delta^{-1} h^2(s_1)/\tilde{A} - \Delta^{-1} h^2(s_o)/\tilde{A} \,\,.
 $$
Again, average this expression, but now over $\Omega(s_1)$. We get
$$ \tilde{G}(s,s_o) \tilde{A} - 0 - G(s,s_o) \tilde{A} + \Delta^{-1}h^2(s) = \int_S \, \Delta^{-1} h^2(s_1) h^2(s_1) \Omega(s_1)/\tilde{A} - \Delta^{-1} h^2(s_o) \,\,.
$$
This gives the transformation formulas $G \rightarrow \tilde{G}$.  The transformation formula  $R \rightarrow \tilde{R}$
follows by the   limit $s \rightarrow s_o$ in
$$ \tilde{G}(s,s_o) - \frac{\log(\tilde{d} (s,s_o))}{2\pi} = \left(\tilde{G}(s,s_o) - G(s,s_o) \right)  +
\left(G(s,s_o) - \frac{\log d(s,s_o)}{2\pi} \right) - \frac{\log(\tilde{d}(s,s_o)/d(s,s_o))}{2\pi} \,\,.
$$

\subsection{Jordan domains, Schottky doubles} 

The ``naive rule'' of Proposition \ref{naive} works in the case of a Jordan domains $D \subset S$, regardless of the vorticity sum.  The basic ingredients for vortex motion  are the {\it hydrodynamical Green functions}:
\begin{definition} A hydrodynamical Green function
 	$G_D(s;s_o)$ is a real harmonic function of $p \in D - s_o$, with (-) logarithmical singularity
at $s = s_o$, {\it symmetric} in $s,s_o$, continuous in $\overline{D} - s_o$.   $G$ is
 constant on the boundaries and with desired circulations around them (topology imposes relations).
 \end{definition}


\noindent C.C. Lin's result on the added term in the Hamiltonian under a conformal transformation is easily interpreted
 in terms of Lemma \ref{naive}, the ``naive rule''.  Let $f:(D,S,g) \rightarrow (\tilde{D},\tilde{S},\tilde{g})$ a conformal map between two Jordan domains,
\begin{equation}
 \tilde{g} ( df(s) \cdot u , df(s) \cdot v ) = h^2(s)\, g(u,v) \,\,.
 \end{equation}
  \begin{proposition} \label{naiverule}
Let $G_D(s;s_o)$ be the hydrodynamical Green function for $D$ . Then:
 \begin{equation} \tilde{G}_{\tilde{D}}(\tilde{s},\tilde{s}_o) =
 G_D(s;s_o)
 \end{equation}
 The regularized Green function receives the naive correction,
 \begin{equation} \tilde{g}(s_o) = g(s_o) - \frac{\log(h)}{2\pi}  \,\,.
 \end{equation}
 \end{proposition}
 Note that the sign in equation (4.6) of \cite{Hally} must be corrected from + to -.


For domains conformal to a planar region  Green functions are classified into two types.  {\it First type}, more employed in electrostatics,  means that $G = 0$ on all boundaries, possibly with different
circulations around each of them. {\it Second kind}, or {\it modified} Green functions are more commonly used
in hydrodynamics.   One prescribes zero circulation around the boundaries,
except at one of them; the values of $G_D = {\rm const.} $  on the boundaries may differ, one of them being normalized to 0.
 The existence of such Green functions For multiply connected planar domains is due to Koebe, and was used by C.C.Lin
in his seminal vortex paper \cite{Lin}.

Specially exciting is the fact that recently Crowdy and Marshall
(\cite{CM},\cite{CM1},\cite{CM2}, \cite{CM3})
 brought to ``implementational fruition'' the task if obtaining Green functions in canonical domains.
 We believe hydrodynamical Green functions should exist on arbitrary Jordan domains on Riemann surfaces - the difference from planar regions is that in general the domain
 will  have ``handles'' (like a juice jar).
 \begin{remark}
An interesting viewpoint is to consider the Schottky double  $S_D = D + \overline{D}$   of a Jordan domain $D$,
$S$ is called a symmetric Riemann surface.  Consider a metric on $S$
such that the reflection is an isometry.  If at   initial time each vortex on $D$ corresponds to a symmetric one in $\overline{D}$, with opposite vorticity, then the corresponding vortices remain symmetric for all time.    N-vortex motion on $D$ is the restriction of the corresponding $2N$ vortex motion on $S$ (using the corresponding Hamiltonian from I).  Thus Schottky double construction   is
reminiscent of Thomson's image vortex method. Inspite of the fact that a given metric in $D$ in general
will not extend to an reflection-isometric one, since the total vorticity is zero,  only the Riemann surface
structure is needed.  In fact,  Schottky-Klein prime functions have been proved useful
to  finding hydrodynamical Green functions on multiply connected planar domains  (\cite{CM11}). 
\end{remark}

\section{When the total vorticity vanishes}

Recall that when the sum of vorticities does not vanish,  using  Green's function for the Laplace-Beltrami operator
on a closed surface    is  unavoidable,  because   there is
a constant {\it nonzero} background counteracting vorticity in $S - \{s_1,...,s_n \}$, and this complicates
the conformal transformation rules.

\subsection{Riemann surface Green functions}  

The case where
$\sum \kappa_j = 0\,$   has a special feature: {\it the collective vortex stream function is the same
for all conformal metrics}.    Of course, this is the basis of conformal mappping methods in applied
hydrodynamics.  Mathematically, this follows from the special   two-dimensional property
$$  \tilde{g} = h^2 g \,\,\,\, \Rightarrow    \Delta_{\tilde{g}} = h^{-2} \Delta_g \,\,.
$$
Hence, irrespective of the chosen metric in its conformal class, all Laplace-Beltrami operators annihilate the same functions,  logarithmic singularities allowed. We need  only   the complex structure ${\cal A}$:  residues at poles are
conformal invariants, etc.

Lset $\, G_S = G(s ;  s_o,s_1) $ be   the {\it two point Green function} for (S,{\cal A}),
whose existence for closed Riemann surfaces $S$ can be used as the starting point\footnote{It is possible that Riemann \cite{Riemann} and Klein
\cite{Klein} insights came from attended Helmoltz lectures.} for the theory of meromorphic functions and differentials  and   uniformization theorems,
 see e.g.  Weyl (\cite{Weyl},II.13, potential arising from a doublet source) or \cite{Springer}.
$G$ is the unique (up to a constant) real {\it harmonic} function for  $s \in S$, except for a (+) logarithmical singularity at $s_o$ and a  (-) logarithmical singularity  at $s_1$.

 $\sum_{i=1}^N \kappa_i = 0$ implies that the collective
vortex stream function can be given in a {\it conformally invariant} way:
\begin{equation}
\psi(s)  =  \psi_{{\rm conformal}}(s) = \sum_{j=1}^N \, \kappa_j \, G(s; s_j, s^*)
\end{equation}
where $s^*$ is any chosen point in $S$.  Irrespective of the metric, the background vorticity of every
term in   {\it zero}  except at the poles $s_j, s^*$. 
Moreover, the choice of $s^*$ is irrelevant as the infinite vorticities there
cancel out because   $\sum_{i=1}^N \kappa_i = 0$.

[Unfortunatelly this is not always a free lunch, since there are some   drawbacks.  One  is that   Riemann surface two
point Green functions have no symmetries in its arguments, which makes it impossible to use $G(s; s_j, s^*)$
for a collective Hamiltonian formulation. Secondly,   there is no direct interpretation
for the desingularization procedure  in terms of invariant objects such as  Robin's function. One can also add a practical objection, that it seems to be equally hard finding (in practice)
$G_{S,g}(s,s_o)$ or  $G_S(s;s_o,s_1)$.  In fact, $G_S(s;s_o,s_1) = G_{S,g}(s,s_o) - G_{S,g}(s,s_1) + {\rm const.}$  where in the
right hand side any metric can be used.]

  At any rate, for any vortex problem in which the background vorticity is zero, the ``naive conformal rule'' (Proposition \ref{naive}) holds, and this fact may simplify matters considerably.

\subsection{Surfaces conformal to $S^2$}  

Take the unit sphere $S^2 \subset \Re^3$ as the ``concrete model'' for $S$ (thought of as an abstract Riemann surface). Any metric can be written as
$g = h^2 g_o \,\,,$
where $g_o$ is the canonical metric on the sphere. Stereographic projection $F: S^2 \rightarrow \C$ from the north pole $s^* = (0,0,1)$ to the equatorial $z$-plane yields \cite{Kimura}
\begin{equation}  g_o = h_o^2 |dz|^2 \,\,, \,\,\,\,  h_o = \frac{2}{1 + |z|^2}  .
\end{equation}
We can use  the   $z$-plane to study the dynamics. If a particle passes through the north pole, in $\C$ it will disappear in the infinite
but reappear instantaneously from  another direction, but ee may assume that this   is not happening during
the time frame we are working with\footnote{{\it Esse est percipi}, stealing the motto of Bishop Berkeley, ironically the greatest foe of Calculus lovers.}.  In $z$-plane, the Green function is $\log(z)/2\pi$.  Let $h \cdot h_o $ the combined conformal
factor from $(S^2, g)$ to $\C$.
Because of zero total vorticity, vortex dynamics on $(S^2,g)$ can be studied in $\C$ with the symplectic form
\begin{equation}
\Omega = \sum_j \,  \kappa_j  h(x_j,y_j)^2 h_o(x_j,y_j)^2\, dx_j \wedge dy_j \,\,\,,
\end{equation}
and Hamiltonian given by Theorem \ref{conformalrule}  (we will omit the factor $2\pi$):
\begin{equation}
H = \frac{1}{4} \sum'_{j,\ell} \kappa_j \kappa_{\ell} \log ( |z_j - z_{\ell}|^2 ) - \sum_{p=1}^N  \, \frac{1}{2} \kappa_{m}^2   \log(h_o(x_m,y_m) h(x_m,y_m))     \,\,\,.
\end{equation}
Here   $\,\,$ ' $\,\,$ means that the $j = \ell$ terms are omitted in the summation.
Now we rewrite each  $\kappa_m^2$ as $- \kappa_m \left( \sum_{p, p \neq m} \kappa_p \right)$, expand and regroup terms. We also use the following interesting formula:
\begin{lemma}  [\cite{Castilho}]
\begin{equation}
|s_1 - s_2|^2 = h_o(s_1)h_o(s_2) |z_1-z_2|^2 \,\, .
\end{equation}
Here the distances are the euclidian distances in $\Re^2$ and $\Re^3$.
\end{lemma}
Pulling back to the sphere, we get:
\begin{proposition} \label{sphere}
For a system of $N$ point vortices in the sphere $S^2$  with a metric $g = h^2 g_o $  ($g_o$ is the standard
metric), and such that {\it the total vorticity vanishes}, the dynamics is governed by
\begin{equation}
\Omega = \sum_j  \,\kappa_j  h^2(s_i) \omega(s_i)\,\,\,,\,\,\,\,
H = \frac{1}{4} \sum'_{j,\ell} \kappa_j \kappa_{\ell}  \log\left( h(s_j)h(s_{\ell}) |s_j - s_{\ell}|^2 \right)
\end{equation}
where  $|s_j - s_{\ell}|$ is the euclidian distance and $\omega$ is the area form of the sphere.
\end{proposition}

\subsection{Proof of Kimura's conjecture on dipole motion}  \label{Kimura}

The Hamiltonian for a vortex pair with opposite vorticities can be written as
\begin{equation} \label{dipoleham}
H = - \kappa^2 \frac{\log d(s_1,s_2)}{2\pi} + \kappa^2  B(s_1,s_2), \,\,\, B(s_1,s_2) = \left[ \frac{R(s_1) + R(s_2)}{2} - \left( G(s_1,s_2) -
\frac{\log d(s_1,s_2)}{2\pi} \right) \right] \,\, .
\end{equation}
 Let
  $\kappa =  O(\epsilon)$ and initial conditions $d(s_1(0),s_2(0)) = O(\epsilon)$.   Kimura's conjectured  that as $\epsilon \rightarrow 0$ the vorticity
center moves on a geodesic.
The truncated system taking only the first term in (\ref{dipoleham}) yields a system of ODEs on $D(s_o) \times D(s_o)$, the product of two copies of a   geodesic ball centered on $s_o \in S$,  given by
 \begin{equation}  \label{vortexpair}
 \dot{s}_1 = - \kappa  \, {\rm sgrad}_{s_1} \, \log d(s_1,s_2)\,\,,\,\, \dot{s}_2 =   \kappa
 {\rm sgrad}_{s_2} \, \log d(s_1,s_2)  \,\,.
 \end{equation}
 We claim that this equation represents the dominant $O(1)$ geodesic dynamics plus an $O(\epsilon^2)$ perturbation. Due to its symmetric form,
  $ B(s_1,s_2)$  also produces
 an $O(\epsilon^2)$ perturbation.
We take Gauss coordinates (\cite{Struik}, \cite{doCarmo})
\begin{equation}  \label{Gausscoordinates}
 ds^2 = du^2 + G(u,v)dv^2\,\,, \,\, G(0,v) = 1\,\,\,,\,\,\, \frac{\partial}{\partial u}_{|u=0} G(u,v) = 0 \,\,.
\end{equation}
The $u$-curves form a field of geodesics, and the central $v$-curve (for $u=0$) is also a geodesic.
Suppose $s_o$ corresponds to some $v$ and  $u=0$.  In these coordinates,  $s_1(0) = (  - \epsilon, v )$
and $s_2(0) = ( \epsilon, v ).$  Clearly, at $t=0$,$\,\,\dot{s}_{1,2}$ will be tangent to $v$-curves, and in the Gauss coordinates,
\begin{equation}
\dot{v}_1(0) = \kappa/2\epsilon \frac{1}{\sqrt{G(-\epsilon, v )}}\,\,\,,\,\,\,
\dot{v}_2(0) = \kappa/2\epsilon \frac{1}{\sqrt{G(\epsilon, v)}} \,\,,\,\,  \dot{u}_1(0) = \dot{u}_2(0) = 0 \,\,.
\end{equation}
In the limit as $\epsilon \rightarrow 0$ (\ref{Gausscoordinates}), with $\kappa = 2\kappa_o \epsilon$,
$\,\, \kappa_o = 1$,
we get
\begin{equation}
\dot{v}_1(0) =  \dot{v}_2(0) = 1 + O(\epsilon)\,\,,\,\,  \dot{u}_1(0), \dot{u}_2(0) = O(\epsilon) \,\,.
\end{equation}
Note that $\frac{\partial}{\partial u}_{|u=0} G(u,v) = 0$ implies that (\ref{vortexpair}) does not
have an $O(\epsilon)$ term.  This concludes the proof.

  \section{Final Remarks}

\subsection{Robin's function and vortex drift}

A single vortex $\dot{s}_o = {\rm sgrad} R(s_o)\,\,$  describes contour lines of Robin's function, so it would be interesting
to understand its generic properties.  For metrics on the sphere $R$ satisfies (see also (\ref{Steiner},  \cite{Steiner})
\begin{equation}
\Delta R = h^2 ( K_E - 4\pi/A(E) )
\end{equation}
where $\Delta$ is the standard Laplacian on the sphere.

  \subsection{Batman's function and integrable  vortex pair problems on Liouville surfaces} 

If  $S$ embedds in $\R^3$ having
an axis of  symmetry $n$, the momentum map of the $S^1$ action on $S \times  ... \times S$ with
 $ \Omega = \sum \kappa_j \Omega_{s_j}$ is  $J = (\sum \kappa_j s_j) \cdot n$.  In particular, a vortex pair on a surface of revolution is integrable, with second integral
$ f = (s_1 - s_2) \cdot n$.
Can we find nontrivial examples of integrable vortex pairs?  As it is well known, Jacobi has shown that the geodesic flow on the ellipsoid  is
integrable. In view of Kimura's theorem,  vortex pairs extends geodesic problems, so the following question is natural:
\begin{question}
Is the vortex pair problem on the triaxial ellipsoid (with opposite vorticities) integrable? If it is not, can one prove it analytically
 by Melnikov's method applied to the geodesics passing by the umbilics?
 \end{question}
Conformal maps from the triaxial ellipsoid  to the standard sphere  are  known
(\cite{Craig}, \cite{Muller}, \cite{Schering}). The conformal factor $h^2$  can be retrieved from these papers, so  Proposition \ref{sphere}
can be implemented.
 Alternatively,  using sphero-conical coordinates (see, eg. \cite{KoillerKurt})
allows obtaining $h^2$ from scratch with two elliptic integrals. More generally, Liouville surfaces (see e.g., \cite{Fomenko}) are those for which the geodesic
flow is integrable.  For which Liouville surfaces is the vortex pair system integrable? In this programme we have found some preliminary results
\cite{Koilleretal} for (\ref{dipoleham}) showing how the geodesic system is perturbed for
$d(s_1,s_2) = O(\epsilon)$.  The $O(\epsilon^2)$ perturbation term  is
$$
 B(s_1,s_2) = \left[ \frac{R(s_1) + R(s_2)}{2} - \left( G(s_1,s_2) -
\frac{\log d(s_1,s_2)}{2\pi} \right) \right]
$$
which could be called, for the sake of justice, Batman's function.


  \subsection{Time dependent problems}
  Either on closed surfaces or Jordan domains it is natural to consider time dependent problems,
  subject
to the natural constraint (in view of imcompressibility) that  the total area  with respect to the metrics
$g(t)$ does not change.
For the sphere and disks Riemann mapping theorem guarantees that the change is always
conformal.  In the general case  the complex structure may also change, an interesting complication that
should be addressed.

 \subsection{Higher dimensional extensions}

 Alan Weinstein called
 our attention that the statement of Theorem \ref{main} makes sense on a compact Kahler manifold,  replacing the two dimensional objects  by their higher dimensional analogues.  In particular, Calabi-Yau manifolds are very fashionable objects nowadays.
 Would the generalization of (\ref{bigsymplectic}, \ref{bighamiltonian})    represent singular solutions to Euler's equation for the hydrodynamics on $(S,g)$, i.e. the geodesic flow on ${\rm SDiff}(S)$, in the spirit of \cite{Arnold} and \cite{MW}? Or will that construction be more in line with
field theorists  approaches, see e.g., Baez \cite{Baez}?

 \subsection{Vortices for a proof of Riemann hypothesis?}

One of the approaches   for a proof of
Riemann hypothesis is  searching a connection between
the notrivial zeros of Riemann's zeta function with some quantum mechanical problem (Berry and Keating, \cite{Berry}).
The quantum version of geodesic motion on $(S,g)$, i.e.,  the free particle, leads to the spectrum of minus the Laplace-Beltrami operator. As a vortex dipole moves on geodesics, why not
 look for a   proof  quantizing a vortex problem, which is
asymmetric with respect to time reversal?

 ``Sir Michael Berry proposes that there exists a classical dynamical system, asymmetric with respect to time reversal, the lengths of whose periodic orbits correspond to the rational primes, and whose quantum-mechanical analog has a Hamiltonian with zeros equal to the imaginary parts of the nontrivial zeros of the zeta function.
The search for such a dynamical system is one approach to proving the Riemann hypothesis   (Daniel Bump, \cite{Bump}).''

 We conclude this note on  a
highly speculative tone:

\begin{conjecture}
A connection may exist between the zeros of  Riemann's zeta function and the  quantization of a 3/2 degrees of freedom vortex monopole problem
on some compact Riemann surface with metric $g$,
$$   \dot{s} = {\rm sgrad} \,G(s,s_o)\,\,\, , \,\,\, \dot{s}_o = {\rm sgrad} R(s_o)\,\,,
$$
or with a vortex pair
$$  \dot{s_1} = - \kappa {\rm sgrad}_{s_1} H(s_1,s_2)\,\,\,\,\, \dot{s_1} =  \kappa {\rm sgrad}_{s_2} H(s_1,s_2)\,\, \,\,\,,\,\,\,
H  = - \kappa^2 \frac{\log d(s_1,s_2)}{2\pi} + \kappa^2  B(s_1,s_2) \,\,\, .
$$
\end{conjecture}

\bigskip
\smallskip

{\bf Acknowlegments.}  The authors thank Alan Weinstein for encouragement and comments, and Jean Singer and Kate Okikiolu for informations about the Laplace-Beltrami Green functions.





\begin{thebibliography}{99}

\bibitem{Aref} Aref, H.,  Point vortex dynamics: A classical mathematics playground, J. Math.  Phys. {\bf 48}:6, 065401  (2007)
\bibitem{Aref1}  Aref, K.,  P. K. Newton, P.K., Stremler, M.A.,  Tokieda, T.,
Vainchtein,  D.L. Vortex Crystals, Adv.   Appl. Math., {\bf 39}, 1--79  (2003)
\bibitem{Arnold}  Arnold, V.I., Sur la geom\'etrie differentielle des groupes de Lie de dimension infinie et ses applications a
 l'hydrodynamique des fluids parfaits,
 Ann. Inst. Grenoble, {\bf 16},  319--361 (1966) 
 \bibitem{Baez}  Baez, J., {\em Knots and Quantum Gravity},  Oxford University Press,  1994
\bibitem{Berry} Berry, M. V.,   Keating, J. P., The Riemann zeros and eigenvalue asymptotics,   SIAM Review {\bf 41},
 236--266 (1999)
 \bibitem{BoattoCrowdy}  Boatto, S., Crowdy, D., Point-vortex dynamics, in
{\em Encyclopedia of Mathematical Physics}, ed. by Irina Arefeva, Daniel Sternheimer, Springer-Verlag, to appear.
\bibitem{BM} Borisov A.V., Mamaev I.S., {\em Mathematical methods in the dynamics of vortex structures}, Moscow - Izhevsk: Institute of Computer Science, 2005   (in Russian)
\bibitem{Bogomolov} Bogomolov, V.A., Dynamics of vorticity at a sphere,
Fluid Dyn. {\bf 6}, 863 -- 870 (1977)
\bibitem{Fomenko}  Bolsinov, A.V.,  Fomenko, A.T., {\em Integrable Hamiltonian Systems}, CRC,  2004
\bibitem{BoPa} Borisov A.V.,  Pavlov, A.E., Dynamics and statics of vortices on a plane and a sphere — I, Reg. Chaotic Dyn., {\bf 3}:1, 28--38 (1998)
\bibitem{BoLe} Borisov A.V.,  Lebedev, V.G., Dynamics of three vortices on a plane and a sphere — II. General compact case, Reg. Chaotic Dyn., {\bf 3}:2,  99--114 (1998)
\bibitem{BoKi} Borisov, A.V., Kilin, A.A., Stability of Thomson's configurations of vortices on a sphere,  Reg. Chaotic Dyn., {\bf 5}:2, 189--100 (2000)
\bibitem{Borisovmass}
 Borisov, A.V.,  Mamaev, I.S.,   Ramodanov, S.M., Dynamics of two interacting circular
cylinders in perfect fluid, Discr. and Cont. Dyn. Systems, {\bf 19}:2, 235-253 (2007)
\bibitem{Bump} Bump, D.,  in http://math.stanford.edu/~bump ;\\ see also
  http://www.secamlocal.ex.ac.uk/people/staff/mrwatkin/zeta/physics1.htm
\bibitem{CaBo} Cabral, H., Boatto, S.,
Nonlinear Stability of a Latitudinal Ring of Point-Vortices on a Nonrotating Sphere, SIAM
J. Appl. Math, {\bf 64}:1, 216--230 (2003)
\bibitem{Castilho} Castilho, C., Machado, H.,  The N-vortex problem on a symmetric ellipsoid: a perturbation approach,
arXiv:[mathDS] 07010.4004v1 (2007)
\bibitem{Craig} Craig, T., Orthomorphic Projection of an Ellipsoid Upon a Sphere,
Amer. J.  Math. {\bf 3}:2, 114-127  (1880)'
\bibitem{CM}  Crowdy, D.,
Point vortex motion on the surface of a sphere with impenetrable boundaries,
Phys. Fluids {18}, 036602, (2006).
\bibitem{CM1}  Crowdy, D., Marshall, J.,
Analytical formulae for the Kirchhoff–Routh path function in multiply connected domains
Proc. R. Soc. A  {\bf 461}, 2477–2501 (2005)
\bibitem{CM11} Crowdy, D.,  Marshall, J.,
Conformal mappings between canonical multiply connected domains,
Comput. Methods Funct. Theory, {\bf 6}:1, 59--76 (2006)
\bibitem{CM2} Crowdy, D., Marshall, J., Green's functions for Laplace's equation in multiply connected domains
 IMA J. Appl. Mathematics  {\bf 72}: 278-301 (2007)
 \bibitem{CM3}  Crowdy, D.,  Marshall, J.,
Computing the Schottky-Klein prime function on the Schottky double of planar domains,
Comput. Methods Funct. Theory, {\bf 7}:1, 293--308 (2007)
\bibitem{doCarmo} do Carmo, M., {\em Differential Geometry of Curves and Surfaces},  Prentice Hall, 1976
 \bibitem{Dong}  Dong, S.,  Kircher, S., Garland,M., Harmonic Functions for Quadrilateral
Remeshing of Arbitrary Manifolds, Computer Aided Geom. Design, {\bf 22}:5, 392--423 (2005)
\bibitem{Euler}  Euler, L. , Principes g\'en\'eraux du mouvement des fluides,  M\'emoires de l'acad\'emie des sciences de Berlin {\bf 11}, pp. 274-315 (1757)
downloadable from  http://math.dartmouth.edu/~euler/pages/E226.html
or from  http://bibliothek.bbaw.de/bibliothek/digital/struktur/02-hist/1755/jpg-0600/00000282.htm  .
\bibitem{Euler2}  Euler, L., Decouverte d'un nouveau principe de Mecanique,
M\'emoires de l'acad\'emie des sciences de Berlin {\bf 6}, 185-217 (1752)
E177 at Euler archives.
\bibitem{Gustafsson} Flucher, M., Gustafsson, B., Vortex motion in 2-dimensional hydrodynamics, energy renormalization and stability of vortex pair, TRITA preprint series (1997)
\bibitem{Hally}  Hally, D., Stability of streets of vortices on surfaces of revolution with a reflection symmetry, J. Math. Phys. {\bf 21}:1, 211--217 (1980)
\bibitem{Lacomba1}
 Hern\`andez-Gardu\~no, A.,  Lacomba, E., Collisions and regularization for the 3-vortex problem,
arXiv:math-ph/0412024 (2004);
Collisions of Four Point Vortices in the Plane,   math-ph/0609016 (2006)
\bibitem{Lacomba2}
 Hern\`andez-Gardu\~no, A.,  Lacomba, E.,
Collisions of Four Point Vortices in the Plane,   math-ph/0609016 (2006)
\bibitem{Helmholtz} Helmoltz, H., \"Uber integrale der hydrodynamischen gleichungen
   welche den Wirbelbewegungen entsprechen,
 Crelles J. {\bf 55}, 25-55 (1858)  downloadable from http://dz-srv1.sub.uni-goettingen.de/sub/digbib/loader?did=D268537
\bibitem{Joe}  Keller, J., Teapot Effect
J. Appl. Phys. {\bf 28:8},  859--864 (1957)
\bibitem{KiNa} Kidambi, R.,   Newton, P.,  Motion of three point vortices on a sphere, Physica D {\bf 116}, 143–175   (1998)
\bibitem{KiN} Kidambi, R., Newton, P.K., Point vortex motion on a sphere with solid boundaries, Phys. Fluids {\bf 12}:3, 581--588 (2000)
\bibitem{Kimura} Kimura, Y., Vortex motion on surfaces with constant curvature, Proc. R. Soc. Lond. A {\bf 455}, 245--259 (1999)
\bibitem{KiO} Kimura, Y., Okamoto, H., Vortex motion on a sphere, J. Phys. Soc. Jps. {\bf 56}, 4203--4206 (1987)
\bibitem{Kirchhoff} Kirchhoff, G., Vorlesungen \"uber mathematische Physik, {\it
Mechanik}, ch. XX, Teubner, 1876.  downloadable from http://gallica.bnf.fr/ark:/12148/bpt6k99608d  .
\bibitem{Klein}  Klein, F., {\em \"Uber Riemann's Theorie der
Algebraischen Functionen}, Gutenberg E-book  20313, 2007
\bibitem{KoillerKurt} Koiller, J., Ehlers, K., Rubber rolling over a sphere, Reg. Chaotic Dyn. {\bf 12}:2, 127--152 (2007)
\bibitem{KRO} Koiller, J., Ragazzo, C., Oliva, W., On the motion of two-dimensional vortices with mass,
J.Nonlinear Science, {\bf 4}:1,  375--418 (1994)
\bibitem{Koilleretal} Koiller, J., Okikiolu, K., Boatto, S., in preparation.
\bibitem {Laurent} Laurent-Polz, F., Point vortices on a rotating sphere,
Reg. Chaotic Dyn., {\bf 10}:1, 39--58 (2005)
\bibitem{Kuchemann}  K\"uchemann, D., Report on the I.U.T.A.M. symposium on concentrated vortex motions in fluids,
 J.   Fluid Mechanics  {\bf 21} 1-20 (1965)
\bibitem{Lim} Lim, C., Montaldi, J., Roberts, M.,
Relative equilibria of point vortices on the sphere,
 Physica D {\bf 148}:1-2,  97--135 (2001)
 \bibitem{Lin} Lin, C.C., On the motion of vortices in two dimensions - I. Existence of the Kirchhoff-Routh function; II - some further investigations on the Kirchhoff-Routh function, Proc. Nat. Ac. Sci. , {\bf 27}, 570-577 (1941)
\bibitem{MW}   Marsden, J.E.,   Weinstein, A.,  Coadjoint orbits, vortices,
 and Clebsch variables for incompressible fluids, Physica D {\bf 7}, 305--323 (1983)   
 \bibitem{Muller} Muller, B., Kartenprojektionen des dreiachsigen Ellipsoids, Diplomarbeit, Geodatisches Institut,
 Univ. Stuttgart (1991)
 \bibitem{Neishtadt}  Neishtadt, A. I., Vasiliev, A. A.,  Destruction of adiabatic invariance at resonances
in slow fast Hamiltonian systems,
 Nuclear Instruments and Methods in Physics Research Section A, {\bf 561}:2, 158--165 (2006)
 \bibitem{Newton} Newton, P., {\em The N-vortex problem. Analytical technques},
Springer, 2001.
\bibitem{Okikiolu} Okikiolu, K., A negative mass theorem for the 2-torus, ArXiv [math-SP]: 07113489v1 (2007)
\bibitem{Okikiolu2} Okikiolu, K., Extremals for logarithmic Hardy-Littlewood-Sobolev inequalities on compact
manifolds, ArXiv [math-SP]: 0603717v2 (2007)
\bibitem{PeMa} Pekarsky, S., Marsden. J., Point vortices on a sphere: stability of relative equilibria, J. Math. Phys. {\bf 39}, 5894 - 5907 (1998)
\bibitem{Polthier}  Polthier, K.,  Preuss,  Identifying Vector Fields Singularities using a
Discrete Hodge Decomposition, in
{\em Visualization and Mathematics III}, Eds: H.C. Hege, K. Polthier, Springer Verlag (2002).
\bibitem{Ramodanov}  Ramodanov, S.M., On the motion of two mass vortices in perfect fluid, pp. 459–-468, in
A.V. Borisov et al. (eds.), {\em IUTAM Symposium on Hamiltonian Dynamics,
Vortex Structures, Turbulence},   Springer,  2007
\bibitem{Riemann} Riemann, B., {Theorie der Abel'schen Functionen}, Journal f\"ur die reine und angewandte Mathematik,
{\bf 54}, 101-155 (1857)
\bibitem{Schering} Schering, E.,  \"Uber die conforme abbildung des ellipsoids auf der
ebene, {\em Gesammelte Mathematische Werke}, ch. III, Mayer and Muller,  Berlin, 1902
\bibitem{Souliere} Souli\`ere, A., Tokieda, T., Periodic motions of vortices on surfaces with symmetry, J. Fluid Mech. {\bf 460}, 83--92 (2002)
\bibitem{Springer}  Springer, G., {\em Introduction to Riemann surfaces}, Addison-Wesley, 1957.
\bibitem{Steiner} Steiner, J., A geometrical mass and its extremal properties for metrics on $S^2$, Duke Math. J.,
{\bf 129}:1, 63--86 (2005)
\bibitem{Struik} Struik, D., {\em Lectures on classical differential geometry}, Addison-Wesley, 1950
\bibitem{Tronin} Tronin, K.G., Absolute choreographies of point vortices on a sphere, Reg. Chaotic Dyn.,  {\bf 11}:1, 123--130 (2006)
\bibitem{Weyl} Weyl, H.,  {\em The concept of a Riemann surface}, Addison-Wesley, 1955



\end{thebibliography}
\end{document}